\documentclass[11pt, reqno]{amsart}

\usepackage{tikz}
\usetikzlibrary{shapes.geometric}

\usepackage{extarrows}

\usepackage{amsxtra}
\usepackage{hyperref}
\usepackage{tikz-cd}
\usepackage[all]{xy}
\usetikzlibrary{patterns}
\usetikzlibrary{backgrounds}

\let\oldtocsection=\tocsection

\let\oldtocsubsection=\tocsubsection

\let\oldtocsubsubsection=\tocsubsubsection

\renewcommand{\tocsection}[2]{\hspace{0em}\oldtocsection{#1}{#2}}
\renewcommand{\tocsubsection}[2]{\hspace{2em}\oldtocsubsection{#1}{#2}}
\renewcommand{\tocsubsubsection}[2]{\hspace{4.5em}\oldtocsubsubsection{#1}{#2}}

\makeatletter
\def\subsection{\@startsection{subsection}{2}
 \z@{.5\linespacing\@plus.7\linespacing}{-.5em}
 {\normalfont\bfseries}}
\def\subsubsection{\@startsection{subsubsection}{3}
 \z@{.5\linespacing\@plus.7\linespacing}{-.5em}
 {\normalfont\bfseries}}
\makeatother

\theoremstyle{plain}

\usepackage{multirow}

\newtheorem{theorem}{Theorem}[section]
\newtheorem{lemma}[theorem]{Lemma}
\newtheorem{definition-theorem}[theorem]{Definition-Theorem}
\newtheorem{proposition}[theorem]{Proposition}
\newtheorem{corollary}[theorem]{Corollary}
\newtheorem{definition}[theorem]{Definition}
\newtheorem{example}[theorem]{Example}
\newtheorem{remark}[theorem]{Remark}
\newtheorem{notation}[theorem]{Notation}
\newtheorem{assumption}[theorem]{Assumption}
\newtheorem{lemma-definition}[theorem]{Lemma-Definition}
\newtheorem{lemma-notation}[theorem]{Lemma-Notation}
\newtheorem{question}[theorem]{Question}
\newtheorem{remark-definition}[theorem]{Remarks-Definition}
\newtheorem{notation-remark}[theorem]{Notation-Remarks}
\newtheorem{conjecture}[theorem]{Conjecture}
\newcommand \bth[1] { \begin{theorem}\label{t#1} }
\newcommand \ble[1] { \begin{lemma}\label{l#1} }

\newcommand \bpr[1] { \begin{proposition}\label{p#1} }
\newcommand \bco[1] { \begin{corollary}\label{c#1} }
\newcommand \bconj[1] { \begin{conjecture}\label{co#1} }
\newcommand \bde[1] { \begin{definition}\label{d#1}\rm }
\newcommand \bex[1] { \begin{example}\label{e#1}\rm }
\newcommand \bre[1] { \begin{remark}\label{r#1}\rm }
\newcommand \bnota[1] {\begin{notation}\label{n#1}\rm }
\newcommand \bas[1] { \begin{assumption}\label{a#1}\rm }

\newcommand \bqu[1] { \begin{question}\label{q#1}\rm }

\newcommand {\eth} { \end{theorem} }
\newcommand {\ele} { \end{lemma} }

\newcommand {\epr} { \end{proposition} }
\newcommand {\econj} { \end{conjecture} }
\newcommand {\eco} { \end{corollary} }
\newcommand {\ede} { \end{definition} }
\newcommand {\eex} { \end{example} }
\newcommand {\ere} { \end{remark} }
\newcommand {\enota} { \end{notation} }
\newcommand {\eas} {\end{assumption}}

\newcommand {\equ} {\end{question}}
\newcommand \thref[1]{Theorem \ref{t#1}}

\newcommand \prref[1]{Proposition \ref{p#1}}

\newcommand \exref[1]{Example \ref{e#1}}

\def \ZZ {{\mathbb Z}}

\def \QQ {{\mathbb Q}}

\newcommand\preceqdot{\mathrel{\ooalign{$\preceq$\cr
  \hidewidth\raise0.225ex\hbox{$\cdot\mkern0.5mu$}\cr}}}

\newcommand{\beqa}{\begin{eqnarray*}}                     
\newcommand{\eeqa}{\end{eqnarray*}}

\def \tkt{{\begin{xy}(0,1)*+{t}="A",(10,1)*+{t'}="B",\ar@{-}^k"A";"B" \end{xy}}}

\def \tit {{\begin{xy}(0,1)*+{t}="A",(10,1)*+{t'}="B",\ar@{-}^i"A";"B" \end{xy}}}
\def \vkv {\begin{xy}(0,1)*+{v}="A",(10,1)*+{v'}="B",\ar@{-}^k"A";"B" \end{xy}}

\usepackage{graphicx} 

\title{When frieze patterns meet Y-systems: {\bf Y}-frieze patterns}
\author{Antoine de Saint Germain}
\address{
Department of Mathematics   \\
The University of Hong Kong \\
Pokfulam Road               \\
Hong Kong}
\email{adsg96@hku.hk}

\begin{document}
\begin{abstract}
    We give an elementary account of the notion of {\bf Y}-frieze patterns, explain some of their properties, and reveal their connection with Coxeter's frieze patterns. 
\end{abstract}
\maketitle

Our objective in this article is to present a variant of 
(Coxeter's) frieze patterns. The latter, introduced by Coxeter in \cite{Cox}, are elementary objects of arithmetic flavour but with a geometric origin traced back to Gauss's {\it pentagramma mirificum}. Coxeter \cite{Cox} proved their ``surprising" periodicity, and Conway and Coxeter \cite{Con-Cox1} explained their remarkable connection with Catalan numbers. In the 2000s, frieze patterns underwent a revival, as they were seen to constitute an ``elementary shadow" of cluster algebras \cite{CC:hall-algebras} of type $A_n$. Since then, much attention has gone into frieze patterns and their various generalisations, from research (e.g. \cite{ARS:friezes, Keller-Sch:linear-recurrence, SOT:2-frieze}) to survey articles (\cite{Sophie-M:survey}, \cite{Baur:surveyFrieze}), lecture series \cite{Pressland:frieze}, Youtube videos on the Numberphile channel \cite{numberphile1, numberphile2}; they even made {\it two} appearances in the {\it Intelligencer} \cite{Tabach:4Vertex, Baur:FriezeIntegers}.

The basic idea of this variant, which we call {\bf Y}-frieze patterns, is similar: an array of numbers satisfying a relation for each ``diamond" of adjacent numbers. The motivation for these patterns comes from 
Zamolodchikov's Y-systems of type $A_n$ (see e.g. \cite{FZ:associahedra, Zamolo:Y-system}). The connection with frieze patterns is a consequence of the general connection between Y-systems and cluster algebras \cite{FG:ensembles,FZ:ClusterIV}; see \cite{Antoine:YFrieze} for more details.

In \S \ref{s:YFrieze}, we define {\bf Y}-frieze patterns, explain how to {\it knit} them and show that they satisfy the same ``surprising" periodicity as frieze patterns. In \S \ref{s:aritPat}, we show that the number of (arithmetic) {\bf Y}-frieze patterns is finite, and formulate the {\it enumeration problem}.  
Finally, in \S \ref{s:fpConnection}, we describe the connection between {\bf Y}-frieze patterns and frieze patterns (whose definition will be recalled), and present a conjecture relevant to the {\it enumeration problem}. 

\section{{\bf Y}-frieze patterns}\label{s:YFrieze}
\subsection{First definitions}

A {\it {\bf Y}-frieze pattern} is a collection of staggered infinite rows of rational numbers arranged as
\[
 \begin{matrix}
    \cdots && 0 && 0 && 0 && 0 && 0 && \cdots \\
    &\cdots &&  b_{0,2} && b_{1,3} && b_{2,4} &&  b_{3,5} && \cdots \\
    \cdots && b_{-1,2}&& b_{0,3} && b_{1,4} && b_{2,5} && b_{3,6} && \cdots\\
    &\cdots &&  b_{-1,3} && b_{0,4} && b_{1,5} &&  b_{2,6} && \cdots \\
    &&&& \vdots && \vdots && \vdots && \\
    
 \end{matrix}
\]
subject to the so-called {\it {\bf Y}-diamond rule:}
\begin{equation}\label{eq:Ydiamond}
 \; WE = (1+N)(1+S)  \quad \text{ for every diamond  } \quad
\begin{matrix}
    & N & \\
    W & & E \\
    & S & 
\end{matrix}.
\end{equation}
The initial row of a {\bf Y}-frieze pattern is called the {\it zeroth} row, and subsequent rows are called the first row, second row, etc. 

\begin{example}
 There is a {\bf Y}-frieze pattern whose every entry in the
 $n^{\rm th}$ row is equal to $(n+1)^2-1$. To see this, note that a typical diamond is of the form \[
    \begin{matrix}
        & n(n-2) & \\
        n^2-1 & & n^2 - 1 \\ 
        & n(n+2) &
    \end{matrix}
\]
and the condition \eqref{eq:Ydiamond} follows from the identity $(n^2-1)(n^2-1) = (n+1)^2(n-1)^2$. The first few rows of this {\bf Y}-frieze pattern are given in Figure \ref{fig:Yfr1}.
\end{example}

\begin{figure}
    \centering
   \[
\begin{matrix}
    && 0 && 0 && 0 && 0 &&  \\
    & 3 && 3 && 3 && 3 \\
    && 8 && 8 && 8 && 8 \\
    & 15 && 15 && 15 && 15 \\
    && \vdots && \vdots && \vdots
\end{matrix}
\] 
\caption{The first four rows of a {\bf Y}-frieze pattern}
    \label{fig:Yfr1}
\end{figure}

By definition, constructing a {\bf Y}-frieze pattern requires choosing (with some care) an infinite array of rationals. We'll now describe two ways in which we can reduce the number of choices.

\subsection{Vertical knitting}
Consider a diamond of the form \eqref{eq:Ydiamond} in which $N, E$ and $W$ are given rationals, and $N \neq -1$. The {\bf Y}-diamond rule implies that 
\begin{equation}\label{eq:vert-knit}
    S = \frac{WE-N-1}{1+N}.
\end{equation}
This observation can be turned into a recursive process called {\it vertical knitting}:

1) choose rational numbers for the first row;

2) knit the second row from the zeroth and first rows using \eqref{eq:vert-knit}; 

3) if the first row does not contain any $-1$s, knit the third row from the first and second rows using \eqref{eq:vert-knit};

4) repeat this process, knitting the $(n+1)^{\rm th}$-row from the $(n-1)^{\rm th}$ and $n^{\rm th}$ rows {\it provided} that the $(n-1)^{\rm th}$ row contains no $-1$s.  

For instance, one can recover the {\bf Y}-frieze pattern in Figure \ref{fig:Yfr1} by vertical knitting the infinite sequence $(3,3,3,3, \ldots)$.

Note that not every {\bf Y}-frieze pattern can be recovered using vertical knitting, due to the possible appearance of $-1$s (see e.g. Figure \ref{fig:Yfri-1s}). 

\begin{figure}[ht]
    \centering
\[
\begin{matrix}
    0&& 0 && 0 && 0 && 0&& 0&& 0 && 0 &&   \\
   & 1 && 1 && 3  && -3&& 0 && 1&& -5 && 1 &&  \\
     &&0 && 2 && -10 && -1 && -1 && -6 && -6 && 0 &&\cdots \\
    &\cdots&& -1 && -6 && -6 && 0 && 2 && -10&& -1 && -1 &&   \\
    &&&& 1 && -5 && 1 && 1 && 3 && -3 && 0 && 1  \\
    &&&&& 0 && 0 && 0 && 0 && 0 && 0 && 0 && 0\\
\end{matrix} 
\]
    \caption{The first six rows of a {\bf Y}-frieze pattern containing $-1$s}
    \label{fig:Yfri-1s}
\end{figure}
In the most extreme case, vertical knitting may produce an entire row of $-1$s, as the following example illustrates. 
\bex{:Yfri-ord5}
Applying the vertical knitting to the sequence $(\, 1,2,5, \,  1,2,5, \, \ldots)$ determines 5 rows. 
\[
\begin{matrix}
    0&& 0 && 0 && 0  \\
   &1&& 2&& 5&& 1 \\
 &&1&& 9&& 4 && 1 \\
 &&&2&& 5&& 1 && 2  && \cdots\\
   &&\cdots&& 0 && 0 && 0 && 0 \\
   &&&&& -1 && -1 && -1 && -1 \\
   &&&&&& 0 && 0 && 0 && 0  \\
   &&&&&&& ? && ? && ? && ?
\end{matrix} 
\]
\eex

If a {\bf Y}-frieze pattern admits a row of $0$s followed by a row of $-1s$, it is said to be {\it closed}. In a closed pattern, if the row of $0$s followed by the row of $-1$s arises in the $(n+1)^{\rm th}$ and $(n+2)^{\rm th}$ rows (for some $n \geq 0$), {\it and} the first, second, $\ldots$ , $n^{\rm th}$ rows are non-zero, this pattern is said to be of {\it width $n$}. When talking about {\bf Y}-frieze patterns of width $n$, we will only consider the rows up to and including the $(n+1)^{\rm th}$ row (the second row of $0$s). For instance, the first five rows of the pattern in \exref{:Yfri-ord5} constitute a closed {\bf Y}-frieze pattern of width 3. For the remainder of this article, we will focus on closed {\bf Y}-frieze patterns.

\subsection{Horizontal knitting} 
We now provide a process for constructing closed {\bf Y}-frieze patterns. Consider once again a diamond of the form \eqref{eq:Ydiamond} in which $N, S$ and $W$ are given {\it positive }rationals. The {\bf Y}-diamond rule implies that 
\begin{equation}\label{eq:hor-knit}
    E = \frac{(1+N)(1+S)}{W}.
\end{equation}
This observation can be turned into a process called {\it horizontal knitting}:

1) place two infinite rows of $0$s, separated by $n$ empty rows (to be filled in by the algorithm); 

2) choose $n$-many {\it positive} rational numbers, and place them on a vertical ``zig-zag" connecting the two rows of $0$s;  

3) knit horizontally using \eqref{eq:hor-knit}, one entry at a time.

\begin{figure}[ht]
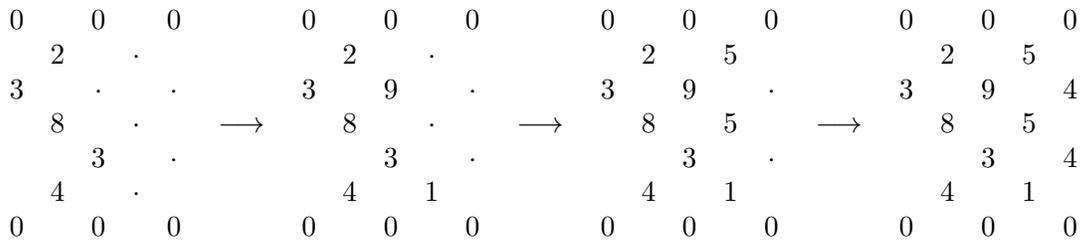

    \centering
    \[
    \begin{matrix}
        0 && 0 && 0 \\
        &2 && \cdot  \\
        3 &&\cdot && \cdot \\
        &8 && \cdot  \\
        &&3 && \cdot  \\
        &4 && \cdot  \\
        0 && 0 && 0 \\
    \end{matrix} 
   \quad \longrightarrow \quad
    \begin{matrix}
        0 && 0 && 0 \\
        &2 && \cdot  \\
        3 &&9 && \cdot \\
        &8 && \cdot  \\
        &&3 && \cdot  \\
        &4 && 1  \\
        0 && 0 && 0 \\
    \end{matrix} 
     \quad \longrightarrow \quad
    \begin{matrix}
        0 && 0 && 0 \\
        &2 && 5  \\
        3 &&9 && \cdot \\
        &8 && 5  \\
        &&3 && \cdot  \\
        &4 && 1  \\
        0 && 0 && 0 \\
    \end{matrix} 
       \quad \longrightarrow \quad
    \begin{matrix}
        0 && 0 && 0 \\
        &2 && 5  \\
        3 &&9 && 4 \\
        &8 && 5  \\
        &&3 && 4  \\
        &4 && 1  \\
        0 && 0 && 0 \\
    \end{matrix} 
\]
    \caption{Horizontal knitting of a {\bf Y}-frieze pattern of width 5}
    \label{fig:horizKnit}
\end{figure}

Here, by vertical ``zig-zag" we mean any choice of one entry per row such that an entry in the $i^{\rm th}$ row is immediately to the south-west or to the south-east of the entry in the $(i-1)^{\rm th}$ row. A few steps of vertical knitting for a given zig-zag are shown in Figure \ref{fig:horizKnit}. 

The following lemma is straightforward.

\ble{}
Fix $n \geq 1$ and a choice of zig-zag. The process of horizontal knitting establishes a bijection 
from $n$-tuples of positive rational numbers to $\QQ_{>0}$-valued {\bf Y}-frieze patterns of width $n$.
\ele

\subsection{Properties of {\bf Y}-frieze patterns}
Our first statement is an analog of Coxeter's celebrated {\it glide symmetry of frieze patterns} (c.f. \S \ref{s:fpConnection} for more on this). 
\bth{:glide}\cite[Theorem 5.2]{Antoine:YFrieze}
    Suppose that $(b_{i,j})$ is a closed {\bf Y}-frieze pattern of width $n$. Then $(b_{i,j})$ is invariant under the {\it glide symmetry} $b_{i,j} = b_{j, i+n+3}$. 
\eth

\begin{figure}[ht]
    \centering
    \includegraphics[scale=0.8]{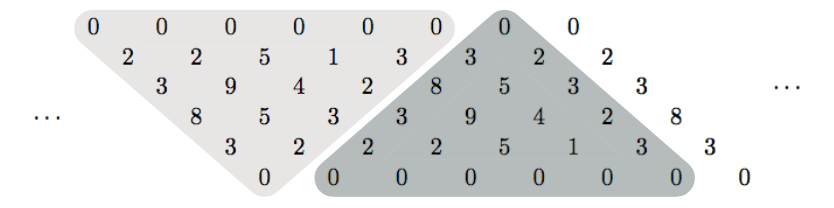}
    \caption{Fundamental domain of a {\bf Y}-frieze pattern of width 4 (light grey) and its image under a glide reflection (dark grey)}
    \label{fig:fundaDom}
\end{figure}

The name ``glide symmetry" is due to the fact that the action of this symmetry on a fundamental domain is by {\it glide reflection}, i.e. composition of a horizontal reflection and a translation. An example of a fundamental domain, the triangle, is illustrated in Figure \ref{fig:fundaDom}. 

Applying the glide symmetry twice, we have that $b_{i,j} = b_{i+n+3, j+n+3}$, i.e. every {\bf Y}-frieze pattern of width $n$ is invariant under (horizontal) translation by $n+3$ steps (see Figure \ref{fig:fundaDom} for an illustration).

\section{Arithmetic patterns}\label{s:aritPat}
Figure \ref{fig:fundaDom} is an example of a {\bf Y}-frieze pattern whose non-zero rows consist entirely of positive integers. We call such patterns {\it arithmetic}. Denote by YFrieze$(n)$ the set of all arithmetic {\bf Y}-frieze patterns of width $n$.

\begin{figure}[ht]
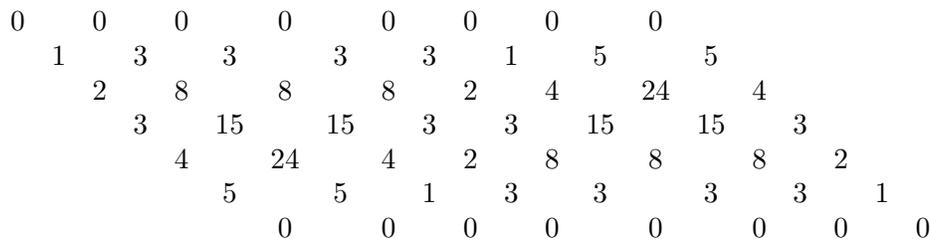

    \centering
    \[
\begin{matrix}
    0 && 0 && 0 && 0 && 0 && 0 && 0 && 0  \\
    & 1 && 3 && 3 && 3 && 3 && 1 && 5 && 5 \\
    && 2 && 8 && 8 && 8 && 2 && 4 && 24 && 4 \\
    &&& 3 && 15 && 15 && 3 && 3 && 15 && 15 && 3 \\
    &&&& 4 && 24 && 4 && 2 && 8 && 8 && 8 && 2  \\
    &&&&&5 && 5 && 1 && 3 && 3 && 3 && 3 && 1 \\
    &&&&&&0 && 0 && 0 && 0 && 0 && 0 && 0 && 0 
\end{matrix}
\]
    \caption{An arithmetic {\bf Y}-frieze pattern of width 5}
    \label{fig:arith-width4}
\end{figure}

\bpr{:unitary}
    The {\bf Y}-frieze pattern of width $n$ obtained by placing the tuple $(1,2, \ldots , n) \in (\ZZ_{>0})^n$ along the north-west to south-east diagonal and performing horizontal knitting, is arithmetic.
\epr
In particular, arithmetic {\bf Y}-frieze patterns of width $n$ exist. An example is given in Figure \ref{fig:arith-width4}.  

\bre{:noLaurent}
Readers familiar with frieze patterns might be tempted to ask whether an $n$-tuple of $1$s can be used to construct arithmetic patterns, instead of the $n$-tuple used in \prref{:unitary}: the resulting {\bf Y}-frieze pattern is in general {\it not} arithmetic, as can be see in Figure \ref{fig:Laurent-phenomenon}. This is because entries in the {\bf Y}-frieze pattern are in general {\it not} Laurent polynomials in the entries of the north-west to south-east diagonal. 
\ere
\begin{figure}[ht]
    \centering
    \[
    \begin{matrix}
        &&0 && 0 && 0 && 0 && 0 && 0\\
        &&& 1 && 2 && \frac{7}{2} && 2 && 1\\
        \cdots&&&& 1 && 6 && 6 && 1 &&&& \cdots\\
        &&&&& 1 && 7 && 1\\
        &&&&&& 0 && 0  
    \end{matrix}
    \]
    \caption{A non-arithmetic {\bf Y}-frieze pattern of width 3 knitted from $1$s on a zig-zag}
    \label{fig:Laurent-phenomenon}
\end{figure}

\prref{:unitary} leads to the very natural challenging question:

\begin{center}
    {\it Enumeration problem: How many arithmetic {\bf Y}-frieze patterns of width $n$ are there? }
\end{center}

A first step towards solving the enumeration problem for {\bf Y}-friezes is given by the following theorem. 

\bth{:finite}\cite[Theorem 5.5]{Antoine:YFrieze}
    Fix $n\geq 1$. The number of arithmetic {\bf Y}-frieze patterns of width $n$ is finite.
\eth

In the next section, we provide a conjectural upper bound on the number of {\bf Y}-frieze patterns for each $n$. For now, let us conclude this section with a table of arithmetic {\bf Y}-frieze patterns of 
width $1, 2$ and $3$. For convenience and readability, we only include the north-west to south-east diagonal of each {\bf Y}-frieze pattern. An interested reader may reconstruct the corresponding patterns by hand, or by using \cite{visualca}.

\begin{center}
\def\arraystretch{1.5}
\begin{tabular}{|c||c|}
\hline
    $n$ & Diagonal of {\bf Y}-frieze patterns of width $n$ \\
    \hline \hline
    1 & $(1)$ \\
    \hline
    2 & $(1,1), (1,2), (2,1), (2,3), (3,2)$  \\
    \hline
    \multirow{2}{*}{3} &  $(1,1,2), (1,2,3), (1,4,5), (2,1,1), (2,3,2)$, \\
    & $(2,9,5), (3,2,1), (3,8,3), (5,4,1), (5,9,2)$.\\
    \hline
\end{tabular}
\end{center}
For $n=1, 2$, these are all {\bf Y}-frieze patterns (\cite[Theorem 4.8]{Antoine:YFrieze}). 
For $n=3$, we conjecture these are all {\bf Y}-frieze patterns.

\section{Connection with frieze patterns}\label{s:fpConnection}
Recall from \cite{Baur:FriezeIntegers} that a {\it frieze pattern} is a collection of staggered infinite rows arranged as
\[
 \begin{matrix}
    \cdots && 1 && 1 && 1 && 1 && 1 && \cdots \\
    &\cdots &&  a_{0,2} && a_{1,3} && a_{2,4} &&  a_{3,5} && \cdots \\
    \cdots && a_{-1,2}&& a_{0,3} && a_{1,4} && a_{2,5} && a_{3,6} && \cdots\\
    &\cdots &&  a_{-1,3} && a_{0,4} && a_{1,5} &&  a_{2,6} && \cdots \\
    &&&& \vdots && \vdots && \vdots && \\
 \end{matrix}
\]
satisfying the {\it diamond rule} $a_{i,j}\, a_{i+1,j+1} = 1 + a_{i,j+1}\, a_{i+1,j}$ for all $i,j$. A frieze pattern is said to be {\it closed} if it has a row of ones followed by a row of zeros, in which case the number $n\geq 1$ of rows strictly between the first two rows of ones is called its {\it width}. As with {\bf Y}-frieze patterns, a closed frieze pattern is said to be {\it arithmetic} if its rows consist entirely of positive integers. Denote by Frieze($n$) the set of all arithmetic frieze patterns of width $n$. It is well-known (see e.g. \cite[Theorem 4.3]{Sophie-M:survey}) that such frieze patterns are counted by Catalan numbers; more precisely we have
\[
\mid\text{Frieze}(n) \mid \, = C_{n+1},
\]
 where $C_{n} = \frac{1}{n+1} \binom{2n}{n}$ is the {\it $n^{\rm th}$ Catalan number}. For example, the number of arithmetic frieze patterns of width $1, 2$ and $3$ is $C_2= 2, C_3 =  5$ and $C_4 = 14$ respectively.  

\bth{:ensemble-map}\cite[Theorem 2.8]{Antoine:YFrieze}
    Fix $n \geq 1$. There is a well-defined map
    \[
        p_n: \text{Frieze}(n) \longrightarrow \text{YFrieze}(n), \qquad (a_{i,j}) \mapsto (b_{i,j}),
    \]
    where $(b_{i,j})$ is the arithmetic {\bf Y}-frieze pattern of width $n$ whose first row is $(a_{i,i+3})_{i \in \ZZ}$, the second row of $(a_{i,j})$. 
\eth
For instance, we have
\[
\begin{matrix}
    1 && 1 && 1 && 1 && 1 &&  \\
    &  2 && 1 && 4 && 1 && 3 &&   \\
    && 1 && 3 && 3 && 2 && 2 &&  \\
    &&& 2 && 2 && 5 && 1 && 3 &&  \\
    &&&& 1 && 3 && 2 &&1 &&4   \\
    &&&&& 1 && 1 && 1 && 1 && 1  
\end{matrix} \;
\stackrel{p_4}{\longmapsto} \; \begin{matrix}
    0 && 0 && 0 && 0 && 0 \\
    & 1 && 3 && 3 && 2 && 2 \\
    && 2 && 8 && 5 && 3 && 3 \\
    &&& 3 && 9 && 4 && 2 && 8 \\
    &&&& 2 && 5 && 1 && 3 && 3 \\
    &&&&& 0 && 0 && 0 && 0 && 0
\end{matrix}
\]

\bre{}
      \thref{:ensemble-map} highlights the importance of the $2^{\rm nd}$ row of a frieze pattern. Note that when $n$ is even, it was observed in \cite[\S 1.8]{SMG:arith-2-frieze} that the $2^{\rm nd}$ row of a frieze pattern of width $n$ can be interpreted as cross-ratios of a configuration of points in $\mathbb{P}^1$.
\ere

The map $p_n$ is the frieze pattern analog of Fock and Goncharov's {\it cluster ensemble map} (also called the {\it p-map}, see \cite{FG:ensembles}). Motivated by properties 
of the cluster ensemble map and an exhaustive search for {\bf Y}-friezes of width up to $12$, we suggest the following conjecture. 
\bconj{:surjective}
    For all $n \geq 1$, the map $p_n$ is surjective.
\econj

Two frieze patterns $f = (f_{i,j})$ and $g = (g_{i,j})$ of width $n$ are said to be {\it {\bf Y}-equivalent}, denoted $f \sim_{\bf Y} g$, if they have the same second row, i.e. if $f_{i,i+3} = g_{i,i+3}$ for all $i \in \ZZ$. Denote by Frieze(n)/$\sim_{\bf Y}$ the set of equivalence classes of $\sim_{\bf Y}$ on Frieze($n$). One can show that 
\begin{itemize}
    \item when $n$ is even, each equivalence class contains a single element (this follows for instance from \cite[Remark 1.18]{SMG:arith-2-frieze}); and 
    \item when $n$ is odd, each equivalence class contains {\it at most} two elements (this follows from \cite[Lemma 7.5]{CuntzHolm:frpatComplex}).
\end{itemize}

Assuming that Conjecture \ref{co:surjective} holds, the discussion above guarantees that 

1) the map $p_n$ descends to a bijection 
\[
    \tilde{p}_n : \text{Frieze}(n)/ \sim_{\bf Y} \quad \longrightarrow  \quad \text{YFrieze}(n),
\]

2) when $n$ is even, $p_n$ is bijective and $|\text{YFrieze}(n)| = C_{n+1}$; and

3) when $n$ is odd, $|\text{YFrieze}(n)| < C_{n+1}$. 

In particular, if Conjecture \ref{co:surjective} holds, the set Frieze(n)/$\sim_{\bf Y}$ is a combinatorial model for arithmetic {\bf Y}-frieze patterns of width $n$. The details of this model, along with an exact formula for the size of Frieze(n)/$\sim_{\bf Y}$ will be explained elsewhere. It would be interesting to find a (direct) combinatorial model for {\bf Y}-frieze patterns (i.e. one which does not involve the use of frieze patterns). 


\textbf{Acknowledgment.} The author would like to thank V. Ovsienko and S. Morier-Genoud for encouragements to write this article. The author has been partially supported by
the Research Grants Council of the Hong Kong SAR, China (GRF 17307718).

\bibliographystyle{alpha}
\bibliography{ref}

\end{document}